\documentclass[12pt,oneside,reqno]{amsart}
\usepackage{bbding}
\usepackage{txfonts}
\usepackage{bbm}
\usepackage{amsmath}
\usepackage{graphicx}
\usepackage{mathrsfs}
\usepackage{stmaryrd}
\usepackage{amsfonts}
\usepackage{enumerate,amsmath,amssymb}

\numberwithin{equation}{section}

\newcommand{\be}{\begin{eqnarray}}
\newcommand{\ee}{\end{eqnarray}}
\newcommand{\ce}{\begin{eqnarray*}}
\newcommand{\de}{\end{eqnarray*}}
\newtheorem{theorem}{Theorem}[section]
\newtheorem{lemma}[theorem]{Lemma}
\newtheorem{remark}[theorem]{Remark}
\newtheorem{definition}[theorem]{Definition}
\newtheorem{proposition}[theorem]{Proposition}
\newtheorem{Examples}[theorem]{Example}
\newtheorem{corollary}[theorem]{Corollary}

\def\eps{\varepsilon}

\def\osc{\mathrm{osc}}
\def\med{\mathrm{med}}

\def\p{\partial}

\def\[{{\Big[}}
\def\]{{\Big]}}
\def\<{{\langle}}
\def\>{{\rangle}}
\def\({{\Big(}}
\def\){{\Big)}}

\def\bx{{\mathbf{x}}}

\def\dif{{\mathord{{\rm d}}}}

\def\no{\nonumber}
\def\={&\!\!=\!\!&}
\def\bt{\begin{theorem}}
\def\et{\end{theorem}}
\def\bl{\begin{lemma}}
\def\el{\end{lemma}}
\def\br{\begin{remark}}
\def\er{\end{remark}}

\def\bd{\begin{definition}}
\def\ed{\end{definition}}
\def\bp{\begin{proposition}}
\def\ep{\end{proposition}}
\def\bc{\begin{corollary}}
\def\ec{\end{corollary}}
\def\bx{\begin{Examples}}
\def\ex{\end{Examples}}
\def\cA{{\mathcal A}}

\def\cF{{\mathcal F}}

\def\cH{{\mathcal H}}
\def\cI{{\mathcal I}}
\def\cJ{{\mathcal J}}

\def\cL{{\mathcal L}}
\def\cM{{\mathcal M}}

\def\cP{{\mathcal P}}

\def\cR{{\mathcal R}}
\def\cS{{\mathcal S}}
\def\cT{{\mathcal T}}

\def\mE{{\mathbb E}}

\def\mH{{\mathbb H}}

\def\mN{{\mathbb N}}

\def\mR{{\mathbb R}}

\def\mU{{\mathbb U}}

\def\mW{{\mathbb W}}
\def\mX{{\mathbb X}}
\def\mY{{\mathbb Y}}

\def\sC{{\mathscr C}}

\def\sE{{\mathscr E}}

\def\sR{{\mathscr R}}

\def\geq{\geqslant}
\def\leq{\leqslant}

\def\div{\mathord{{\rm div}}}

\allowdisplaybreaks

\begin{document}

\title{Well-posedness of fully nonlinear and nonlocal critical parabolic equations$^*$}

\date{}
\author{Xicheng Zhang}


\thanks{$*$This work is supported by NSFs of China (No. 10971076).}

\address{
School of Mathematics and Statistics,
Wuhan University, Wuhan, Hubei 430072, P.R.China\\
Email: XichengZhang@gmail.com
 }

\begin{abstract}
In this paper we prove the existence of smooth solutions to fully nonlinear and nonlocal parabolic equations
with critical index. The proof relies on the apriori H\"older estimate for advection fractional-diffusion
equation established by Silvestre \cite{Si2}.
\end{abstract}

\maketitle \rm

\section{Introduction and main result}

In this paper we are interested in solving the following fully nonlinear and nonlocal parabolic equation:
\begin{align}
\p_t u=F(t,x,u,\nabla u,-(-\Delta)^{\frac{\alpha}{2}}u),\ \ u(0)=\varphi,\ \ \alpha\in(0,2),\label{EE1}
\end{align}
where $F(t,x,u,w,q): [0,1]\times\mR^d\times\mR\times\mR^d\times\mR\to\mR$ is a measurable function,
and $(-\Delta)^{\frac{\alpha}{2}}$ is the usual fractional Laplacian defined by
$$
(-\Delta)^{\frac{\alpha}{2}}u=\cF^{-1}(|\cdot|^\alpha\cF u),\ \ u\in\cS(\mR^d),
$$
where $\cF$ denotes the Fourier's transform, $\cS(\mR^d)$ is the Schwartz class of smooth
real-valued rapidly decreasing functions.

Recently, in the sense of viscosity solutions, fully nonlinear and nonlocal elliptic and parabolic equations have
been extensively studied (cf. \cite{Ca-Si,Si3,Ba-Ch-Im,La-Da}, etc.). In \cite{Ca-Si}, Caffarelli and Silvestre studied
the following type of nonlocal equation:
$$
I_\alpha u(x):=\sup_i\inf_j \left(c^{ij}+b^{ij}\cdot\nabla u(x)
+\int_{\mR^d}[u(x+y)-u(x)]a^{ij}(y)|y|^{-d-\alpha}\dif y\right)=0,
$$
where $\alpha\in(0,2)$, $i,j$ ranges in arbitrary sets, $c^{ij}\in\mR$ and $b^{ij}\in\mR^d$,
the kernel $a^{ij}(y)$ satisfies
$$
a^{ij}(y)=a^{ij}(-y),\ \ a_0\leq a^{ij}(y)\leq a_1.
$$
This type of equation appears in the stochastic control problems. In \cite{Ca-Si},
the extremal Pucci operators are used to characterize the ellipticity,
and the ABP estimate, Harnack inequality and interior
$C^{1,\beta}$-regularity were obtained. In \cite{Si2}, Silvestre studied the following nonlocal
parabolic equation with critical index $\alpha=1$:
$$
\p_tu=I_1u,\ \ u(0)=\varphi,
$$
and established $C^{1,\beta}$-regularity of viscosity solutions. In particular, the following first order
Hamilton-Jacobi equation is covered by the above equation when $H$ is Lipschitz continuous:
$$
\p_t u+H(\nabla u)+(-\Delta)^{\frac{1}{2}}u=0.
$$
In \cite{La-Da}, Lara and Davila extended Silvestre's result to the more general case, and in particular,
focused on the uniformity of regularity as $\alpha\to 2$.

However, it is not known how to solve the fully nonlinear and nonlocal equation (\ref{EE1})
in Sobolev spaces. Let us fix the main idea of the present paper for solving (\ref{EE1}).
Assume that $F$ does not depend on $u$. Taking the gradient with respect to $x$ for equation (\ref{EE1}), we have
$$
\p_t\nabla u=-(\p_q F)(-\Delta)^{\frac{\alpha}{2}}\nabla u+(\nabla_wF)\nabla^2u+\nabla_x F.
$$
We make the following observation:
$$
-(-\Delta)^{\frac{\alpha}{2}}u=(-\Delta)^{\frac{\alpha-2}{2}}\div \nabla u=\cR^\alpha\cdot\nabla u,
$$
where $\cR^\alpha:=(-\Delta)^{\frac{\alpha-2}{2}}\div$ is a bounded linear operator from Bessel potential space
$\mH^{\alpha-1,p}$ to $L^p$ provided $p>1$. If we set $w:=\nabla u$, then $w$
satisfies the following quasi-linear parabolic system:
\begin{align}
\p_tw=-(\p_qF)(w,\cR^\alpha w)(-\Delta)^{\frac{\alpha}{2}}w+(\nabla_wF)(w,\cR^\alpha w)\nabla w
+(\nabla_x F)(w,\cR^\alpha w).\label{EQ3}
\end{align}
It is noticed that the classical quasi-geostrophic equation takes the same form (cf. \cite{Con, Ca-Va, Ki-Na-Vo}, etc.):
$$
\p_t\theta+(-\Delta)^{\frac{\alpha}{2}}\theta+\cR\theta\cdot\nabla\theta=0,\ \ \cR:=\nabla^\bot(-\Delta)^{-\frac{1}{2}}.
$$
Assume now that one can solve equation (\ref{EQ3}), then it is natural to define
$$
u(t,x):=\varphi(x)+\int^t_0 F(s,x,w(s,x),\cR^\alpha w(s,x))\dif s.
$$
Thus, if one can show
\begin{align}
\nabla u=w,\label{EQ4}
\end{align}
then it follows that
$$
u(t,x)=\varphi(x)+\int^t_0 F(s,x,\nabla u(s,x),-(-\Delta)^{\frac{\alpha}{2}}u(s,x))\dif s.
$$
For solving equation (\ref{EQ3}), we shall use Silvestre's H\"older estimate \cite{Si2}
about the following linear parabolic equation:
\begin{align}
\p_tu=-a(-\Delta)^{\frac{\alpha}{2}}u+b\cdot\nabla u+f.\label{EH6}
\end{align}
For proving (\ref{EQ4}), we need to solve a linear equation like
\begin{align}
\p_t u=a(-\Delta)^{\frac{\alpha-2}{2}}\square u+b\cdot(\nabla u-(\nabla u)^t),\label{EE4}
\end{align}
where $\square:=\div \nabla-\nabla\div$ is a symmetric operator on $L^2(\mR^d;\mR^d)$ and
$$
\<\square u,u\>_2=-\|\nabla u\|^2_2+\|\div u\|^2_2.
$$
Notice that in one dimensional case, $\square=0$.

In this work, we mainly concentrate on the critical case $\alpha=1$ and prove the following result:
\bt
\label{Main} Assume that $\p_qF\geq a_0>0$ and for some $\kappa_0>0$,
\begin{align}
|F(t,x,u,0,0)|\leq\kappa_0(|u|+1);\label{EW3}
\end{align}
and for any $R>0$,
\begin{align}
F&\in L^\infty([0,1];C^\infty_b(\mR^d\times B^1_R\times B^d_R\times B^1_R)),\label{EW4}\\
\p_q F,\nabla_wF&\in L^\infty([0,1];C^1_b(\mR^d\times B^1_R\times B^d_R\times\mR)),\label{EW5}\\
\p_uF&\in L^\infty([0,1]\times\mR^d\times B^1_R\times\mR^d\times\mR),\label{EW6}
\end{align}
where $B^d_R$ denotes the open ball in $\mR^d$ with radius $R$ and center $0$;
and for any $j\in\mN$ and $R>0$, there exist $C_{R,j}\geq 0$, $\gamma_{R,j}\geq 0$ and
$h_{R,j}\in (L^1\cap L^\infty)(\mR^d)$ such that
for all $(t,x,u,w,q)\in[0,1]\times\mR^d\times B^1_R\times\mR^d\times\mR$,
\begin{align}
|\nabla^j_x F(t,x,u,w,q)|\leq C_{R,j}|w|(|w|^{\gamma_{R,j}}+1)+h_{R,j}(x),\label{EW7}
\end{align}
where $\gamma_{R,1}=0$. Then for any initial value $\varphi\in \mU^\infty:=\cap_{k,p}\mU^{k,p}$,
where $\mU^{k,p}$ is defined by (\ref{EM7}) below,
there exists a unique $u\in C([0,1];\mU^\infty)$ solving equation (\ref{EE1}) with $\alpha=1$.
Moreover,
$$
\sup_{t\in[0,1]}\|u(t)\|_\infty\leq e^{\kappa_0}(\|\varphi\|_\infty+\kappa_0).
$$
\et
\br
Let $A(q)\in C^\infty(\mR)$ have bounded derivatives of first and second orders and
$\p_q A$ be bounded below by $a_0>0$.
Let $H\in C^\infty(\mR^d)$ and $f\in C^\infty(\mR)$ satisfy $|f(u)|\leq\kappa_0(|u|+1)$. Then
$$
F(t,x,u,w,q):=A(q)+H(w)+f(u)
$$
satisfies all the conditions (\ref{EW3})-(\ref{EW7}).
\er

In the subcritical case $\alpha\in(1,2)$, when we adopt the same argument described above to
solve the fully nonlinear equation (\ref{EE1}), there are two difficulties occurring:
on one hand, we need to prove a stronger apriori H\"older estimate for equation (\ref{EH6})
(see Theorem \ref{Th6} below)
$$
\sup_{t\in[0,1]}\sup_{x\not=y}\frac{|u(t,x)-u(t,y)|}{|x-y|^\beta}\leq C,\ \ \exists\beta\in(\alpha-1,1),
$$
where $C$ only depends on the bounds of $a,b, f$ and $u(0)$;
on the other hand, for $\alpha\in(1,2)$, it is not known whether the uniqueness holds
for equation (\ref{EE4}) in the class of smooth solutions. In the case of $\alpha\in(0,1]$,
this problem can be solved by observing $\div\square u=0$ (see Lemma \ref{Le4}).

In the supercritical case $\alpha\in(0,1)$, it is well-known that there exists an explosion solution
for one-dimensional fractal Burger's equation (cf. \cite{Ki-Na-Sh, Si2}). Nevertheless,
from the proof of Theorem \ref{Main}, one can see that the approach also works for
the following fully nonlinear equation:
$$
\p_t u=F(t,x,u,-(-\Delta)^{\frac{\alpha}{2}}u),\ \ u(0)=\varphi.
$$

The paper is organized as follows: In Section 2, we prepare some notations and recall some well-known facts
for later use. In Section 3, we solve the linear equation in Sobolev spaces. In Section 4,
we prove the existence of smooth solutions for the quasi-linear nonlocal parabolic system. In Section 5,
we give the proof of Theorem \ref{Main}.
\section{Preliminaries}

Let $\mN_0:=\mN\cup\{0\}$. For $p\in(1,\infty)$ and $\beta\in\mN_0$,
let $\mW^{\beta,p}$ be the completion of $\cS(\mR^d)$
with respect to the norm
$$
\|f\|_{\beta,p}:=\sum_{k=0}^\beta\|\nabla^k f\|_p,
$$
where $\nabla^k$ denotes the $k$-order gradient;
and for $0<\beta\not=\mbox{integer}$, let $\mW^{\beta,p}$ be the completion of $\cS(\mR^d)$
with respect to the norm
\begin{align}
\|f\|_{\beta,p}:=\|f\|_p+\sum_{k=0}^{[\beta]}
\left(\iint_{\mR^d\times\mR^d}\frac{|\nabla^kf(x)-\nabla^kf(y)|^p}{|x-y|^{d+\{\beta\}p}}\dif x\dif y\right)^{\frac{1}{p}},\label{Eq4}
\end{align}
where for a number $\beta>0$, $[\beta]$ denotes the integer part of $\beta$
and $\{\beta\}=\beta-[\beta]$. It is well-known that for $k\leq m$, $\theta\in(0,1)$
with $(1-\theta)k+m\theta\notin\mN$ (cf. \cite[p.185, (2)]{Tr}),
\begin{align}
(\mW^{k,p},\mW^{m,p})_{\theta,p}=\mW^{(1-\theta)k+m\theta,p},\label{EU8}
\end{align}
where $(\cdot,\cdot)_{\theta,p}$ stands for the real interpolation space. For $t\in[0,1]$,
write $\mY^{k,p}_t:=L^p([0,t];\mW^{k,p})$ with the norm
$$
\|u\|_{\mY^{k,p}_t}:=\left(\int^t_0\|u(s)\|^p_{k,p}\dif s\right)^{\frac{1}{p}},
$$
and let $\mX^{k,p}_t$ be the completion of all functions $u\in C^\infty([0,t];\cS(\mR^d))$ with
respect to the norm
$$
\|u\|_{\mX^{k,p}_t}:=\sup_{s\in[0,t]}\|u(s)\|_{k-1,p}+\|u\|_{\mY^{k,p}_t}+\|\p_tu\|_{\mY^{k-1,p}_t}.
$$
It is well-known that (cf. \cite[p.180, Theorem III 4.10.2]{Am})
\begin{align}
\mX^{k,p}_t\hookrightarrow C([0,t];\mW^{k-\frac{1}{p},p}).\label{Em}
\end{align}
Let $\mU^{k,p}$ be the Banach space of the completion of $C^\infty_b(\mR^d)$ with respect to the norm:
\begin{align}
\|f\|_{\mU^{k,p}}:=\|f\|_\infty+\|\nabla f\|_{k,p}.\label{EM7}
\end{align}
For simplicity of notation, we also write
$$
\mX^{k,p}:=\mX^{k,p}_1,\ \ \mY^{k,p}:=\mY^{k,p}_1
$$
and
$$
\mW^\infty:=\cap_{k,p}\mW^{k,p},\ \ \mY^\infty:=\cap_{k,p}\mY^{k,p},\ \
\mX^\infty:=\cap_{k,p}\mX^{k,p},\ \ \mU^\infty:=\cap_{k,p}\mU^{k,p}.
$$

Let $\Omega$ be an open domain of $\mR^d$. For $k\in\mN_0\cup\{\infty\}$,
we use $C^k_b=C^k_b(\Omega)$ to denote the space of all bounded and
$k$-order continuous differentiable functions with all bounded derivatives up to $k$-order.
For $\beta\in(0,1)$, let $\sC^\beta$ be the completion of $\cS(\mR^d)$ with respect to the norm
$$
\|f\|_{\sC^\beta}:=\|f\|_\infty+|f|_{\sC^\beta},
$$
where $\|\cdot\|_\infty$ is the sup-norm and
\begin{align}
|f|_{\sC^\beta}:=\sup_{x\not= y}\frac{|f(x)-f(y)|}{|x-y|^{\beta}}.\label{Hold}
\end{align}
Notice that $C^\infty_b(\mR^d)\not\subset\sC^\beta$.
By the Sobolev embedding theorem, one has
$$
\mW^{1,p}\subset \sC^{1-\frac{d}{p}},\ \ p>d.
$$

Let $\sR$ be the class of all linear operators $\cR: \mW^\infty\to\mW^\infty$ satisfying
that for each $\beta\geq 0$ and $p>1$,
$$
\cR: \mW^{\beta,p}\to\mW^{\beta,p} \mbox{ is a bounded linear operator},
$$
and for each $\beta\in(0,1)$,
\begin{align}
|\cR f|_{\sC^\beta}\leq C_{d,\beta}|f|_{\sC^\beta},\ \ \forall f\in\sC^\beta.\label{Ho}
\end{align}
A typical example of such an operator is the Riesz transform:
$$
\cR_jf:=(-\Delta)^{-\frac{1}{2}}\p_jf=\lim_{\eps\to 0}\int_{|y|>\eps}f(x-y)\frac{y_j}{|y|^{d+1}}\dif y.
$$
Indeed, it holds that for any $p>1$ (cf. \cite{St}),
\begin{align}
\|\nabla f\|_p\simeq\|(-\Delta)^{\frac{1}{2}}f\|_p.\label{Riesz}
\end{align}

Recalling that for any $f\in\cS(\mR^d)$,
\begin{align}
(-\Delta)^{\frac{1}{2}}f(x)=c_d\int_{\mR^d}[f(x)-f(x+y)]|y|^{-d-1}\dif y,\label{EM3}
\end{align}
where $c_d>0$ is a universal constant,
we have
\begin{align}
(-\Delta)^{\frac{1}{2}}(fg)=g(-\Delta)^{\frac{1}{2}}f+f(-\Delta)^{\frac{1}{2}} g-\sE(f,g), \label{For}
\end{align}
where
\begin{align}
\sE(f,g)(x):=c_d\int_{\mR^d}(f(x)-f(x+y))(g(x)-g(x+y))|y|^{-d-1}\dif y.\label{For1}
\end{align}
From this formula, it is easy to derive that (see \cite{Zh2}),
\bl\label{Le1}
Let $\zeta\in \cS(\mR^d)$ and set $\zeta_z(x):=\zeta(x-z)$ for $z\in\mR^d$.
Then for any $p\in[1,\infty)$, there exists a constant $C=C(p,d)>0$ such that for all $ f\in\mW^{1,p}$,
\begin{align}
\int_{\mR^d}\|(-\Delta)^{\frac{1}{2}}( f\zeta_z)
-(-\Delta)^{\frac{1}{2}}f\zeta_z\|_p^p\dif z
\leq C\|\zeta\|_{2,p}^p\|f\|_p^{p/2}\|f\|_{1,p}^{p/2}.\label{EE3}
\end{align}
\el

For given $\lambda_0>0$, $f\in L^\infty([0,1];\mW^\infty)$ and $\varphi\in\mW^\infty$,
let us consider the following heat equation:
$$
\p_t u+\lambda_0(-\Delta)^{\frac{1}{2}}u=f,\ \ u(0)=\varphi.
$$
It is well-known that the unique solution can be represented by
$$
u(t,x)=\cP^{\lambda_0}_t\varphi(x)+\int^t_0\cP^{\lambda_0}_{t-s}f(s,x)\dif s,
$$
where $(\cP^{\lambda_0}_t)_{t\geq 0}$ is the Cauchy semigroup associated with
$\lambda_0(-\Delta)^{\frac{1}{2}}$ and given by
$$
\cP^{\lambda_0}_t \varphi(x):=c_d t\int_{\mR^d}\frac{\varphi(\lambda_0y+x)}{(|y|^2+t^2)^{(d+1)/2}}\dif y,
$$
where $c_d>0$ is a universal constant. By the classical Littlewood-Paley-Stein theory,
there exists a constant $C>0$ only depending on $\lambda_0,d,p$
such that for any $f\in L^p([0,1]\times\mR^d)$ (cf. \cite{St, Zh2}),
\begin{align}
\int^1_0\left\|\nabla\int^t_0\cP^{\lambda_0}_{t-s}f(s)\dif s\right\|_p^p\dif s\leq C\int^1_0\|f(s)\|^p_p\dif s.\label{ET3}
\end{align}

We now use the probabilistic technique to extend the above estimate to the more general case.
Let $(L_t)_{t\geq 0}$ be a $d$-dimensional Cauchy process with L\'evy measure
$\nu(\dif x)=\dif x/|x|^{d+1}$. It is well-known that (cf. \cite{Ap})
$$
\cP^{\lambda_0}_t\varphi(x)=\mE\varphi(x+\lambda_0 L_t).
$$
Let $\vartheta:[0,1]\to\mR^d$ and $\lambda:[0,1]\to[0,\infty)$ be two bounded measurable functions. Define
\begin{align}
\cT^{\lambda,\vartheta}_{t,s}\varphi(x):=\mE \varphi\left(x-\int^t_s\vartheta(r)\dif r+\int^t_s\lambda(r)\dif L_r\right).
\label{EY2}
\end{align}
By the theory of stochastic differential equation (cf. \cite[p.402, Theorem 6.7.4]{Ap}), one knows that
$$
\p_t\cT^{\lambda,\vartheta}_{t,s}\varphi+\lambda(-\Delta)^{\frac{1}{2}}\cT^{\lambda,\vartheta}_{t,s}\varphi
+\vartheta\cdot\nabla\cT^{\lambda,\vartheta}_{t,s}\varphi=0.
$$
Now if we define
$$
u(t,x):=\int^t_0\cT^{\lambda,\vartheta}_{t,s}f(s,x)\dif s,
$$
then it is easy to see that
$$
\p_t u+\lambda(-\Delta)^{\frac{1}{2}} u+\vartheta\cdot\nabla u=f,\ \ u(0)=0.
$$
We have
\bt\label{Th4}
Let $\vartheta:[0,1]\to\mR^d$ and $\lambda:[0,1]\to[\lambda_0,\infty)$ be two bounded measurable functions,
where $\lambda_0>0$. For any $p\in(1,\infty)$, there exists a constant
$C$ depending only on $\lambda_0,d,p$ such that for all $f\in L^p([0,1]\times\mR^d)$,
$$
\int^1_0\left\|\nabla\int^t_0\cT^{\lambda,\vartheta}_{t,s}f(s)\dif s\right\|_p^p\dif s\leq C\int^1_0\|f(s)\|^p_p\dif s.
$$
\et
\begin{proof}
Let $(L^{(i)}_t)_{t\geq 0}, i=1,2$ be two independent copies of Cauchy process $(L_t)_{t\geq 0}$.
By the theory of stochastic differential equation (cf. \cite{Ap, Zh4}), one can write
\begin{align}
\cT^{\lambda,\vartheta}_{t,s}\varphi(x)&=\mE \varphi\left(x-\int^t_s\vartheta(r)\dif r
+\int^t_s(\lambda(r)-\lambda_0)\dif L^{(1)}_r+\lambda_0(L^{(2)}_t-L^{(2)}_s)\right)
=\mE\cP^{\lambda_0}_{t-s}\varphi\left(x-X_t+X_s\right),\label{EY1}
\end{align}
where $\cP^{\lambda_0}_t\varphi(x):=\mE\varphi(x+\lambda_0L^{(2)}_t)$ is
the semigroup associated with $\lambda_0(-\Delta)^{\frac{1}{2}}$, and
$$
X_t:=\int^t_0\vartheta(r)\dif r-\int^t_0(\lambda(r)-\lambda_0)\dif L^{(1)}_r.
$$
Define
$$
u(t,x):=\int^t_0\cP^{\lambda_0}_{t-s}f\left(s,x+X_s\right)\dif s.
$$
By (\ref{EY1}) one has
$$
\int^t_0\cT^{\lambda,\vartheta}_{t,s}f(s,x)\dif s=\mE u\left(t,x-X_t\right).
$$
Hence, by H\"older's inequality and Fubini's theorem,
\begin{align*}
\int^1_0\left\|\nabla\int^t_0\cT^{\lambda,\vartheta}_{t,s}f(s,x)\dif s\right\|_p^p\dif t
&=\int^1_0\big\|\mE \nabla u\left(t,\cdot-X_t\right)\big\|_p^p\dif t
\leq\mE\int^1_0\big\|\nabla u\left(t,\cdot-X_t\right)\big\|_p^p\dif t\\
&=\mE\int^1_0\|u(t)\|^p_p\dif t
=\mE\int^1_0\left\|\nabla \int^t_0\cP^{\lambda_0}_{t-s}f\left(s,\cdot+X_s\right)\dif s\right\|_p^p\dif t\\
&\stackrel{(\ref{ET3})}{\leq} C\mE\int^1_0\big\|f\left(s,\cdot+X_s\right)\big\|_p^p\dif s=C\int^1_0\|f(s)\|^p_p\dif s.
\end{align*}
The proof is finished.
\end{proof}

Below we prove a maximum principle for the fully nonlinear equation (\ref{EE1}).
\bt\label{Max}
(Maximum principle) Let $F(t,x,w,q): [0,1]\times\mR^d\times\mR^d\times\mR\to\mR$ be
a measurable function. Assume that for any $R>0$ and all $(t,x,w,q)\in [0,1]\times\mR^d\times\mR^d\times\mR$
with $|w|,|q|\leq R$,
\begin{align}
0\leq \p_q F(t,x,w,q)\leq a_{R,1},\ \ |\nabla_wF|(t,x,w,q)\leq a_{R,2},\label{EG1}
\end{align}
where $a_{R,1}, a_{R,2}>0$. Let $u\in C([0,1]; C^2_b(\mR^d))$ satisfy
$$
u(t,x)=u(0,x)+\int^t_0 F(s,x,\nabla u(s,x),-(-\Delta)^{\frac{1}{2}}u(s,x))\dif s.
$$
If $F(s,x,0,0)\leq 0$, then for all $t\in[0,1]$,
\begin{align}
\sup_{x\in\mR^d}u(t,x)\leq \sup_{x\in\mR^d}u(0,x).\label{EW9}
\end{align}
In particular,
\begin{align}
\|u(t)\|_\infty\leq\|u(0)\|_\infty+\int^t_0 \|F(s,\cdot,0,0)\|_\infty\dif s.\label{EW8}
\end{align}
\et
\begin{proof}
First of all, we assume
\begin{align}
F(s,x,0,0)\leq\delta<0.\label{EM5}
\end{align}
Suppose that (\ref{EW9}) does not hold, then there exists a time $t_0\in(0,1]$ such that
$$
\sup_{x\in\mR^d}u(t_0,x)=\sup_{t\in[0,1]}\sup_{x\in\mR^d}u(t,x).
$$
Let $x_n\in\mR^d$ be such that
$$
\lim_{n\to\infty}u(t_0,x_n)=\sup_{x\in\mR^d}u(t_0,x)\geq u(t,x),\ \ \forall (t,x)\in[0,1]\times\mR^d.
$$
We have for any $\eps\in(0,t_0)$,
\begin{align}
0&\leq\frac{1}{\eps}\left[\lim_{n\to\infty}u(t_0,x_n)-\varlimsup_{n\to\infty}u(t_0-\eps,x_n)\right]\\
&\leq\frac{1}{\eps}\varliminf_{n\to\infty}(u(t_0,x_n)-u(t_0-\eps,x_n))\no\\
&=\frac{1}{\eps}\varliminf_{n\to\infty}\int^{t_0}_{t_0-\eps}
F(s,x_n,\nabla u(s,x_n),-(-\Delta)^{\frac{1}{2}}u(s,x_n))\dif s,\label{EM1}
\end{align}
and for any $h\in\mR^d$,
\begin{align*}
0&\leq\lim_{\eps\downarrow 0}\frac{1}{\eps}\left(\lim_{n\to\infty}u(t_0,x_n)-\varlimsup_{n\to\infty} u(t_0,x_n-\eps h)\right)\\
&\leq\lim_{\eps\downarrow 0}\frac{1}{\eps}\varliminf_{n\to\infty}(u(t_0,x_n)-u(t_0,x_n-\eps h))\\
&=\lim_{\eps\downarrow 0}\varliminf_{n\to\infty}\int^1_0h\cdot\nabla u(t_0,x_n-\eps sh)\dif s\\
&=\varliminf_{n\to\infty}(h\cdot\nabla u(t_0,x_n)).
\end{align*}
In particular, by the arbitrariness of $h$, we get
\begin{align}
\lim_{n\to\infty}\nabla u(t_0,x_n)=0.\label{EM2}
\end{align}
On the other hand, since for any $y\in\mR^d$,
$$
u(t_0,x_n+y)-u(t_0,x_n)\leq\sup_{x\in\mR^d}u(t_0,x)-u(t_0,x_n)\to 0,
$$
by (\ref{EM3}) we have
\begin{align}
\varlimsup_{n\to\infty}-(-\Delta)^{\frac{1}{2}}u(t_0,x_n)
\leq c_d\int_{\mR^d}\varlimsup_{n\to\infty}[u(t_0,x_n+y)-u(t_0,x_n)]|y|^{-d-1}\dif y\leq 0.\label{EM4}
\end{align}
Moreover, since by $u\in C([0,1]; C^2_b(\mR^d))$,
$$
\lim_{s\to t_0}\|\nabla (u(s)-u(t_0))\|_\infty=0
$$
and
$$
\lim_{s\to t_0}\|(-\Delta)^{\frac{1}{2}}(u(s)-u(t_0))\|_\infty=0,
$$
we have by (\ref{EG1}),
$$
\lim_{\eps\to 0}\frac{1}{\eps}\int^{t_0}_{t_0-\eps}
\|F(s,\nabla u(s),-(-\Delta)^{\frac{1}{2}}u(s))-
F(s,\nabla u(t_0),-(-\Delta)^{\frac{1}{2}}u(t_0))\|_\infty\dif s=0.
$$
Hence, by (\ref{EM1}), (\ref{EM2}) and (\ref{EM5}),
\begin{align}
0&\leq\lim_{\eps\to 0}\frac{1}{\eps}\varliminf_{n\to\infty}\int^{t_0}_{t_0-\eps}
F(s,x_n,\nabla u(t_0,x_n),-(-\Delta)^{\frac{1}{2}}u(t_0,x_n))\dif s\no\\
&=\lim_{\eps\to 0}\frac{1}{\eps}\varliminf_{n\to\infty}\int^{t_0}_{t_0-\eps}
F(s,x_n,0,-(-\Delta)^{\frac{1}{2}}u(t_0,x_n))\dif s\no\\
&\leq\lim_{\eps\to 0}\frac{1}{\eps}\varliminf_{n\to\infty}\int^{t_0}_{t_0-\eps}
[F(s,x_n,0,-(-\Delta)^{\frac{1}{2}}u(t_0,x_n))-F(s,x_n,0,0)]\dif s+\delta\no\\
&=\lim_{\eps\to 0}\varliminf_{n\to\infty}
\Big[a_{n,\eps}\cdot\Big(-(-\Delta)^{\frac{1}{2}}u(t_0,x_n)\Big)\Big]+\delta,\label{EM6}
\end{align}
where
$$
a_{n,\eps}:=\frac{1}{\eps}\int^{t_0}_{t_0-\eps}\int^1_0
(\p_qF)(s,x_n,0,-r(-\Delta)^{\frac{1}{2}}u(t_0,x_n))\dif r\dif s.
$$
Let $R:=\|(-\Delta)^{\frac{1}{2}}u(t_0)\|_\infty$. Noticing that
$$
0\leq a_{n,\eps}\leq a_{R,1},
$$
by (\ref{EM4}), (\ref{EM6}) and $\delta<0$, we obtain a contradiction.

We now drop assumption (\ref{EM5}). For $\delta<0$, set
$$
u_\delta(t,x)=u(t,x)+\delta t.
$$
Then
$$
u_\delta(t,x)=u_\delta(0,x)+\int^t_0 \Big[\delta+F(s,x,\nabla u_\delta(s,x),-(-\Delta)^{\frac{1}{2}}u_\delta(s,x))\Big]\dif s.
$$
So, for all $t\in[0,1]$,
$$
\sup_{x\in\mR^d}u(t,x)\leq \sup_{x\in\mR^d}u_\delta(t,x)-\delta t
\leq\sup_{x\in\mR^d}u(0,x)-\delta t.
$$
Letting $\delta\uparrow 0$, we conclude the proof of (\ref{EW9}).

As for (\ref{EW8}), by considering
$$
\tilde u(t,x)=u(t,x)-\int^t_0\|F(s,\cdot,0,0)\|_\infty\dif s
$$
and using (\ref{EW9}) for $\tilde u(t,x)$ and $-\tilde u(t,x)$ respectively, we immediately obtain (\ref{EW8}).
\end{proof}

Next we recall Silvestre's H\"older estimate about the linear advection fractional-diffusion equation.
The following result is taken
from \cite[Corollary 6.2]{Zh2}.
Although the proofs given in \cite{Si2} and \cite{Zh2} are only for
constant diffusion coefficient $a(t,x)$, by slight modifications, they are also adapted
to the general bounded measurable function $a(t,x)$.
\bt\label{Th6}
(Silvestre's H\"older estimate)  Let $a:[0,1]\times\mR^d\to\mR$ and $b:[0,1]\times\mR^d\to\mR^d$ be two bounded measurable
functions. Let $u\in C([0,1];C^2_b(\mR^d))$ satisfy
$$
u(t,x)=u(0,x)-\int^t_0(a(-\Delta)^{\frac{1}{2}}u)(s,x)\dif s+\int^t_0(b\cdot\nabla u)(s,x)\dif s+\int^t_0 f(s,x)\dif s.
$$
If $a(t,x)\geq a_0>0$, then for any $\gamma\in(0,1)$, there exist
$\beta\in(0,1)$ and $C>0$ depending only on $d, a_0,\gamma$ and $\|a\|_\infty, \|b\|_\infty$ such that
\begin{align}
\sup_{t\in[0,1]}|u(t)|_{\sC^\beta}\leq C(\|u\|_\infty+\|f\|_\infty+|u(0)|_{\sC^\gamma}),\label{EY9}
\end{align}
where $|\cdot|_{\sC^\beta}$ is defined by (\ref{Hold}).
\et

\section{Linear nonlocal parabolic equation}

In this section, we consider the following linear scalar nonlocal equation:
\begin{align}
\p_t u+a(-\Delta)^{\frac{1}{2}} u+b\cdot\nabla u=f,\ \ u(0)=\varphi,\label{Eq1}
\end{align}
where $a:[0,1]\times\mR^d\to\mR$ and $b:[0,1]\times\mR^d\to\mR^d$ are two bounded measurable functions.

An increasing function $\omega:\mR^+\to\mR^+$ is called a modulus function if $\lim_{s\downarrow 0}\omega(s)=0$.
We make the following assumptions on $a$ and $b$:
\begin{enumerate}[{\bf (H$^{a,b}_k$)}]
\item Let $k\in\mN_0$, $a,b\in L^\infty([0,1];C^k_b)$, and
there are two modulus functions $\omega_a$ and $\omega_b$ such that
for all $t\in[0,1]$ and $x,y\in\mR^d$,
\begin{align}
|a(t,x)-a(t,y)|\leq\omega_a(|x-y|),\ \ |b(t,x)-b(t,y)|\leq\omega_b(|x-y|).\label{Eq2}
\end{align}
Moreover, for some $a_0,a_1>0$ and all $(t,x)\in[0,1]\times\mR^d$,
$$
a_0\leq a(t,x)\leq a_1.
$$
\end{enumerate}

We first prove the following important apriori estimate.
\bl\label{Le7}
For given $p\in(1,\infty)$ and $k\in\mN$, let $f\in \mY^{k-1,p}$ and
$u\in\mX^{k,p}$ satisfy that for almost all $(t,x)\in[0,1]\times\mR^d$,
\begin{align}
\p_tu(t,x)+a(t,x)(-\Delta)^{\frac{1}{2}} u(t,x)+b(t,x)\cdot\nabla u(t,x)=f(t,x).\label{EE5}
\end{align}
Then under {\bf (H$^{a,b}_{k-1}$)}, there exists a constant $C_{k,p}>0$ such that for all $t\in[0,1]$,
\begin{align}
\|u\|_{\mX^{k,p}_t}\leq C_{k,p}\left(\|u(0)\|_{k-\frac{1}{p},p}+\|f\|_{\mY^{k-1,p}_t}\right),\label{Es8}
\end{align}
where $C_{1,p}$ depends only on $a_0,a_1,\|b\|_\infty,d,p$ and $\omega_a,\omega_b$.
In particular, equation (\ref{EE5}) admits at most one solution in $\mX^{k,p}$.
\el
\begin{proof}
Let $(\rho_\eps)_{\eps\in(0,1)}$ be a family of mollifiers in $\mR^d$, i.e.,
$\rho_\eps(x)=\eps^{-d}\rho(\eps^{-1}x)$,
where $\rho\in C^\infty_0(\mR^d)$ is nonnegative and has support in $B_1$ and $\int\rho=1$. Define
$$
u_\eps(t):=u(t)*\rho_\eps,\ \ a_\eps(t):=a(t)*\rho_\eps,\ \ b_\eps(t):=b(t)*\rho_\eps,\ \ f_\eps(t):=f(t)*\rho_\eps.
$$
Taking convolutions for both sides of (\ref{EE5}), we have
\begin{align}
\p_tu_\eps(t,x)+a_\eps(t,x)(-\Delta)^{\frac{1}{2}} u_\eps(t,x)+b_\eps(t,x)\cdot\nabla u_\eps(t,x)=h_\eps(t,x),\label{Eq333}
\end{align}
where
\begin{align*}
h_\eps(t,x)&:=f_\eps(t,x)+a_\eps(t,x)(-\Delta)^{\frac{1}{2}} u_\eps(t,x)-
[(a(t)(-\Delta)^{\frac{1}{2}} u(t))*\rho_\eps](x)\\
&\quad+b_\eps(t,x)\cdot\nabla u_\eps(t,x)-[(b(t)\cdot\nabla u(t))*\rho_\eps](x).
\end{align*}
By (\ref{Eq2}), it is easy to see that for all $\eps\in(0,1)$ and $t\in[0,1]$ and $x,y\in\mR^d$,
\begin{align}
|a_\eps(t,x)-a_\eps(t,y)|\leq\omega_a(|x-y|),\ \
|b_\eps(t,x)-b_\eps(t,y)|\leq\omega_b(|x-y|),\label{Eq5}
\end{align}
and
$$
|a_\eps(t,x)-a(t,x)|\leq\omega_a(\eps),\ \ |b_\eps(t,x)-b(t,x)|\leq\omega_b(\eps).
$$
Moreover, by the property of convolutions, we also have
$$
\lim_{\eps\to0}\int^1_0\|h_\eps(t)-f(t)\|^p_p\dif t=0.
$$
Below, we use the method of freezing the coefficients to prove that for all $t\in[0,1]$,
\begin{align}
\|u_\eps\|_{\mX^{1,p}_t}\leq C\left(\|u_\eps(0)\|_{1-\frac{1}{p},p}+\|h_\eps\|_{\mY^{0,p}_t}\right),\label{Ep7}
\end{align}
where the constant $C$ is independent of $\eps$. After proving this estimate, (\ref{Es8}) with $k=1$
immediately follows by taking limits for (\ref{Ep7}).

For simplicity of notation, we drop the subscript $\eps$ below.
Fix $\delta>0$ being small enough, whose value will be determined below.
Let $\zeta$ be a smooth function with support in $B_\delta$ and $\|\zeta\|_p=1$. For $z\in\mR^d$, set
$$
\zeta_z(x):=\zeta(x-z),\ \ \lambda^a_z:=a(t,z),\ \ \vartheta^b_z(t):=b(t,z).
$$
Multiplying both sides of (\ref{Eq333}) by $\zeta_z$, we have
$$
\p_t(u\zeta_z)+\lambda^a_z(-\Delta)^{\frac{1}{2}} (u\zeta_z)+\vartheta^b_z\cdot\nabla(u\zeta_z)=g^\zeta_z,
$$
where
\begin{align*}
g^\zeta_z&:=(\lambda^a_z-a)(-\Delta)^{\frac{1}{2}} u\zeta_z
+\lambda^a_z((-\Delta)^{\frac{1}{2}}(u\zeta_z)-(-\Delta)^{\frac{1}{2}} u\zeta_z)\\
&\quad+(\vartheta^b_z-b)\cdot\nabla (u\zeta_z)+ub\cdot\nabla\zeta_z+h\zeta_z.
\end{align*}
Let $\cT^{\lambda^a_z,\vartheta^b_z}_{t,s}$ be defined by (\ref{EY2}).
Then $u\zeta_z$ can be written as
$$
u\zeta_z(t,x)=\cT^{\lambda^a_z,\vartheta^b_z}_{t,0}(u(0)\zeta_z)(x)
+\int^t_0\cT^{\lambda^a_z,\vartheta^b_z}_{t,s} g^\zeta_z(s,x)\dif s,
$$
and so that for any $T\in[0,1]$,
\begin{align*}
\int^T_0\|\nabla(u\zeta_z)(t)\|^p_p\dif t&\leq
2^{p-1}\left(\int^T_0\|\nabla\cT^{\lambda^a_z,\vartheta^b_z}_{t,0}(u(0)\zeta_z)\|_p^p\dif t
+\int^T_0\left\|\nabla\int^t_0\cT^{\lambda^a_z,\vartheta^b_z}_{t,s}
g^\zeta_z(s)\dif s\right\|^p_p\dif t\right)\\
&=:2^{p-1}(I_1(T,z)+I_2(T,z)).
\end{align*}
For $I_1(T,z)$, by (\ref{EY1}) and (\ref{Riesz}), we have
\begin{align}
\int^T_0\|\nabla\cT^{\lambda^a_z,\vartheta^b_z}_{t,0}(u(0)\zeta_z)\|^p_p\dif t
&=\int^T_0\left\|\nabla\cP^{a_0}_t(u(0)\zeta_z)\right\|^p_p\dif t\no\\
&\leq C\int^T_0\left\|(-\Delta)^{\frac{1}{2}}\cP^{a_0}_t(u(0)\zeta_z)\right\|^p_p\dif t\no\\
&\leq C\|u(0)\zeta_z\|^p_{1-\frac{1}{p},p},\label{ET2}
\end{align}
where the last step is due to \cite[p.96 Theorem 1.14.5]{Tr} and (\ref{EU8}).
Thus, by definition (\ref{Eq4}), it is easy to see that
$$
\int_{\mR^d}I_1(T,z)\dif z\leq C\int_{\mR^d}\|u(0)\zeta_z\|^p_{1-\frac{1}{p},p}\dif z
\leq C\Big(\|u(0)\|^p_{1-\frac{1}{p},p}\|\zeta\|_p^p+\|u(0)\|_p^p\|\zeta\|^p_{1-\frac{1}{p},p}\Big).
$$
For $I_2(T,z)$, by Theorem \ref{Th4}, we have
\begin{align*}
I_2(T,z)&\leq C\int^T_0\|g^\zeta_z(s)\|^p_p\dif s\leq
C\int^T_0\|((\lambda^a_z-a)((-\Delta)^{\frac{1}{2}} u\zeta_z))(s)\|^p_p\dif s\\
&\qquad+C\int^T_0\|\lambda^a_z((-\Delta)^{\frac{1}{2}} u\zeta_z-(-\Delta)^{\frac{1}{2}}(u\zeta_z))(s)\|^p_p\dif s\\
&\qquad+C\int^T_0\|((\vartheta^b_z-b)\cdot\nabla(u\zeta_z))(s)\|^p_p\dif s\\
&\qquad+C\int^T_0\|(ub\cdot\nabla\zeta_z)(s)\|^p_p\dif s+C\int^T_0\|h\zeta_z(s)\|^p_p\dif s\\
&=:I_{21}(T,z)+I_{22}(T,z)+I_{23}(T,z)+I_{24}(T,z)+I_{25}(T,z).
\end{align*}
For $I_{21}(T,z)$, by (\ref{Eq5}) and $\|\zeta\|_p=1$, we have
\begin{align*}
\int_{\mR^d}I_{21}(T,z)\dif z&\leq C\omega^p_a(\delta)\int_{\mR^d}\int^T_0\|((-\Delta)^{\frac{1}{2}} u\zeta_z)(s)\|^p_p\dif s\dif z\\
&= C\omega^p_a(\delta)\int^T_0\|(-\Delta)^{\frac{1}{2}} u(s)\|^p_p\dif s
\stackrel{(\ref{Riesz})}{\leq} C\omega^p_a(\delta)\int^T_0\|\nabla u(s)\|^p_p\dif s.
\end{align*}
For $I_{22}(T,z)$, by (\ref{EE3}) and Young's inequality, we have
\begin{align*}
\int_{\mR^d}I_{22}(T,z)\dif z&\leq
Ca_1\int^T_0\int_{\mR^d}\|((-\Delta)^{\frac{1}{2}} u\zeta_z-(-\Delta)^{\frac{1}{2}}(u\zeta_z))(s)\|^p_p\dif z\dif s\\
&\leq C\int^T_0\|u(s)\|^p_p\dif s+C\int^T_0\|u(s)\|^{p/2}_p\|\nabla u(s)\|^{p/2}_p\dif s\\
&\leq C\int^T_0\|u(s)\|^p_p\dif s+\frac{1}{4^p}\int^T_0\|\nabla u(s)\|^p_p\dif s.
\end{align*}
For $I_{23}(T,z)$, as above we have
$$
\int_{\mR^d}I_{23}(T,z)\dif z\leq C\omega^p_b(\delta)\left(\int^T_0\|\nabla u(s)\|^p_p\dif s+\|\nabla\zeta\|^p_p\int^T_0\|u(s)\|^p_p\dif s\right).
$$
Moreover, it is easy to see that
\begin{align*}
\int_{\mR^d}I_{24}(T,z)\dif z&\leq C\|b\|^p_\infty\|\nabla\zeta\|_p^p\int^T_0\|u(s)\|^p_p\dif s,\\
\int_{\mR^d}I_{25}(T,z)\dif z&\leq C\int^T_0\|h(s)\|^p_p\dif s.
\end{align*}
Combining the above calculations, we get
\begin{align*}
&\int^T_0\|\nabla u(s)\|^p_p\dif s=\int^T_0\int_{\mR^d}\|\nabla u(s)\cdot\zeta_z\|^p_p\dif z\dif s\\
&\quad\leq 2^{p-1}\int^T_0\int_{\mR^d}\|\nabla(u\zeta_z)(s)\|^p_p\dif z\dif s
+2^{p-1}\|\nabla\zeta\|_p^p\int^T_0\|u(s)\|^p_p\dif s\\
&\quad\leq C\|u(0)\|^p_{1-\frac{1}{p},p}+\Big(\frac{1}{4}+C(\omega^p_a(\delta)+\omega^p_b(\delta))\Big)\int^T_0\|\nabla u(s)\|^p_p\dif s\\
&\quad\quad+C\int^T_0\|u(s)\|^p_p\dif s+C\int^T_0\|h(s)\|^p_p\dif s.
\end{align*}
Choosing $\delta_0>0$ being small enough so that
$$
C(\omega^p_a(\delta_0)+\omega^p_b(\delta_0))\leq \frac{1}{4},
$$
we obtain that for all $T\in[0,1]$,
\begin{align}
\int^T_0\|\nabla u(s)\|^p_p\dif s\leq C\|u(0)\|^p_{1-\frac{1}{p},p}+
C\int^T_0\|u(s)\|^p_p\dif s+C\int^T_0\|h(s)\|^p_p\dif s.\label{EL1}
\end{align}
On the other hand, by (\ref{Eq333}), it is easy to see that for all $t\in[0,1]$,
$$
\|u(t)\|^p_p\leq C\|u(0)\|_p^p+C\int^t_0 \|\nabla u(s)\|^p_p\dif s+C\int^t_0 \|h(s)\|^p_p\dif s,
$$
which together with (\ref{EL1}) and Gronwall's inequality yields that for all $t\in[0,1]$,
\begin{align}
\sup_{s\in[0,t]}\|u(s)\|^p_p+\int^t_0\|\nabla u(s)\|^p_p\dif s\leq C
\left(\|u(0)\|^p_{1-\frac{1}{p},p}+\int^t_0 \|h(s)\|^p_p\dif s\right).\label{EI1}
\end{align}
From equation (\ref{EE5}), by (\ref{Riesz}) we also have
\begin{align*}
\int^t_0\|\p_su(s)\|^p_p\dif s&\leq C\left(\|a\|_\infty^p\int^t_0\|(-\Delta)^{\frac{1}{2}} u(s)\|^p_p\dif s
+\|b\|_\infty^p\int^t_0\|\nabla u(s)\|^p_p\dif s
+\int^t_0\|h(s)\|^p_p\dif s\right)\\
&\leq C\left((\|a\|_\infty^p+\|b\|_\infty^p)\int^t_0\|\nabla u(s)\|^p_p\dif s
+\int^t_0\|h(s)\|^p_p\dif s\right),
\end{align*}
which together with (\ref{EI1}) gives (\ref{Ep7}), and therefore (\ref{Es8}) with $k=1$.

Let us now estimate the higher order derivatives. For $n=1,2,\cdots,k$, let
$$
w^{(n)}(t,x):=\nabla^n u(t,x).
$$
By the chain rule, we have
$$
\p_tw^{(n)}+a(-\Delta)^{\frac{1}{2}} w^{(n)}+b\cdot\nabla w^{(n)}=h^{(n)},
$$
where
$$
h^{(n)}:=\nabla^n f-\sum_{j=1}^n \frac{n!}{(n-j)!j!}\Big(\nabla^ja\cdot\nabla^{n-j}(-\Delta)^{\frac{1}{2}} u+
\nabla^jb\cdot\nabla^{n-j+1}u\Big).
$$
Thus, by (\ref{Es8}) with $k=1$ and the assumptions, we have
\begin{align*}
\|\nabla^nu\|_{\mX^{1,p}_t}=\|w^{(n)}\|_{\mX^{1,p}_t}&\leq C\left(\|w^{(n)}(0)\|_{1-\frac{1}{p},p}+\|h^{(n)}\|_{\mY^{0,p}_t}\right)\\
&\leq C\left(\|\nabla^n u(0)\|_{1-\frac{1}{p},p}+\sum_{j=1}^n\|\nabla^{n-j+1}u\|_{\mY^{0,p}_t}
+\|\nabla^nf\|_{\mY^{0,p}_t}\right),
\end{align*}
which implies that
$$
\|u\|_{\mX^{n+1,p}_t}\leq C\left(\|u(0)\|_{n+1-\frac{1}{p},p}+\|u\|_{\mY^{n,p}_t}+\|f\|_{\mY^{n,p}_t}\right).
$$
By induction method, one obtains (\ref{Ep7}).
\end{proof}

Now we prove the existence of solutions to equation (\ref{Eq1}).
\bt\label{Th1}
Let $p>1$ and $k\in\mN$. Under {\bf (H$^{a,b}_{k-1}$)}, for any $\varphi\in\mW^{k-\frac{1}{p},p}$
and $f\in\mY^{k-1,p}$, there exists a unique $u\in \mX^{k,p}$ with $u(0)=\varphi$ solving equation (\ref{Eq1}).
\et
\begin{proof}
We use the continuity method. For $\lambda\in[0,1]$, define an operator
$$
U_\lambda:=\p_t+\lambda a(-\Delta)^{\frac{1}{2}}+\lambda b\cdot\nabla+(1-\lambda)(-\Delta)^{\frac{1}{2}}.
$$
By (\ref{Riesz}), it is easy to see that
\begin{align}
U_\lambda: \mX^{k,p}\to \mY^{k-1,p}.\label{Eq6}
\end{align}
For given $\varphi\in\mW^{k-\frac{1}{p},p}$, let $\mX^{k,p}_\varphi$ be the space
of all functions $u\in\mX^{k,p}$ with $u(0)=\varphi$. It is clear that $\mX^{k,p}_\varphi$
is a complete metric space with respect to the metric $\|\cdot\|_{\mX^{k,p}}$.
For $\lambda=0$ and $f\in \mY^{k-1,p}$, it is well-known that
there is a unique $u\in\mX^{k,p}_\varphi$ such that
$$
U_0 u=\p_t u+(-\Delta)^{\frac{1}{2}} u=f.
$$
In fact, by Duhamel's formula, the unique solution can be represented by
$$
u(t,x)=\cP^1_t\varphi(x)+\int^t_0\cP^1_{t-s}f(s,x)\dif s.
$$
Suppose now that for some $\lambda_0\in[0,1)$, and for any $f\in \mY^{k-1,p}$, the equation
$$
U_{\lambda_0}u=f
$$
admits a unique solution $u\in\mX^{k,p}_\varphi$.
Then, for fixed $f\in\mY^{k-1,p}$ and $\lambda\in[\lambda_0,1]$, and for any $u\in\mX^{k,p}_\varphi$,
by (\ref{Eq6}), the equation
\begin{align}
U_{\lambda_0}w=f+(U_{\lambda_0}-U_\lambda)u\label{Eq0}
\end{align}
admits a unique solution $w\in\mX^{k,p}_\varphi$.
Introduce an operator
$$
w=Q_\lambda u.
$$
We want to use Lemma \ref{Le7} to show that there exists an $\eps>0$
independent of $\lambda_0$ such that for all $\lambda\in[\lambda_0,\lambda_0+\eps]$,
$$
Q_\lambda:\mX^{k,p}_\varphi\to\mX^{k,p}_\varphi
$$
is a contraction operator.

Let $u_1,u_2\in\mX^{k,p}_\varphi$ and $w_i=Q_{\lambda}u_i,i=1,2$. By equation (\ref{Eq0}), we have
\begin{align*}
U_{\lambda_0}(w_1-w_2)=(U_{\lambda_0}-U_\lambda)(u_1-u_2)=(\lambda_0-\lambda)((a-1)(-\Delta)^{\frac{1}{2}}+b\cdot\nabla)(u_1-u_2).
\end{align*}
By (\ref{Es8}) and (\ref{Riesz}), one sees that
\begin{align*}
\|Q_\lambda u_1-Q_\lambda u_2\|_{\mX^{k,p}}
&\leq C_{k,p}|\lambda_0-\lambda|\cdot\|((a-1)(-\Delta)^{\frac{1}{2}}+b\cdot\nabla)(u_1-u_2)\|_{\mY^{k-1,p}}\\
&\leq C_0|\lambda_0-\lambda|\cdot\|u_1-u_2\|_{\mY^{k,p}}
\leq C_0|\lambda_0-\lambda|\cdot\|u_1-u_2\|_{\mX^{k,p}},
\end{align*}
where $C_0$ is independent of $\lambda,\lambda_0$ and $u_1,u_2$. Taking $\eps=1/(2C_0)$, one sees that
$$
Q_\lambda: \mX^{k,p}_\varphi\to \mX^{k,p}_\varphi
$$
is a $1/2$-contraction operator. By the fixed point theorem, for each $\lambda\in[\lambda_0,\lambda_0+\eps]$,
there exists a unique $u\in\mX^{k,p}_\varphi$ such that
$$
Q_\lambda u=u,
$$
which means that
$$
U_\lambda u=f.
$$
Now starting from $\lambda=0$, after repeating the above construction
$[\frac{1}{\eps}]+1$-steps, one obtains that for any $f\in\mY^{k-1,p}$,
$$
U_1u=f
$$
admits a unique solution $u\in\mX^{k,p}_\varphi$.
\end{proof}

\section{Quasi-linear nonlocal parabolic system}

Consider the following quasi-linear nonlocal parabolic system:
\begin{align}
\p_t u+a(u,\cR_a u)(-\Delta)^{\frac{1}{2}} u+ b(u,\cR_b u)\cdot\nabla u=f(u,\cR_f u),\label{EQ1}
\end{align}
where $u=(u^1,\cdots,u^m)$ and
\begin{align*}
&a(t,x,u,r):[0,1]\times\mR^d\times\mR^m\times\mR^k\to\mR,\\
&b(t,x,u,r):[0,1]\times\mR^d\times\mR^m\times\mR^k\to\mR^d,\\
&f(t,x,u,r):[0,1]\times\mR^d\times\mR^m\times\mR^k\to\mR^m,
\end{align*}
are measurable functions, and
\begin{align*}
&\cR_a=(\cR^{ij}_a), \cR_b=(\cR^{ij}_b), \cR_f=(\cR^{ij}_f)\in\sR^{k\times m}.
\end{align*}
Here we have used that $\cR_a u=\sum_{j=1}^m\cR_a^{ij} u^j$, similarly for $\cR_b u$ and $\cR_f u$.

The main result of this section is:
\bt\label{Th2}
Suppose that $a(t,x,u,r)\geq a_0>0$, and for any $R>0$,
\begin{align}
a,b,f&\in L^\infty([0,1];C^\infty_b(\mR^d\times B^m_R\times B_R^k)),\label{EO3}\\
a,b&\in L^\infty([0,1];C^1_b(\mR^d\times B^m_R\times \mR^k)),\label{EO33}
\end{align}
where $B^m_R$ denotes the ball in $\mR^m$ with radius $R$,
and for each $j\in\mN_0$, there exist $C_{f,j},\gamma_j\geq 0$ and $h_j\in (L^1\cap L^\infty)(\mR^d)$ such that
\begin{align}
|\nabla^j_x f(t,x,u,r)|\leq C_{f,j}|u|(|u|^{\gamma_j}+1)+h_j(x),\label{EO1}
\end{align}
and for some $C_f\geq 0$,
\begin{align}
\<u, f(t,x,u,r)\>_{\mR^m}\leq C_f(|u|^2+1).\label{EO2}
\end{align}
Then for any $\varphi\in\mW^\infty$, there exists a unique $u\in\mX^\infty$ solving equation (\ref{EQ1}). Moreover,
$$
\sup_{t\in[0,1]}\|u(t)\|^2_\infty\leq e^{C_f}(\|\varphi\|^2_\infty+C_f).
$$
\et
\begin{proof}
First of all, for any $\cR\in\sR$ and $u\in\mX^{k,p}$,
by the boundedness of $\cR$ in Sobolev space $\mW^{k,p}$, one has
$$
(t,x)\mapsto \cR u(t,x)\in \mX^{k,p}.
$$
Thus, by (\ref{EO3}) and the chain rules, one sees that for any $u\in\mX^\infty$,
\begin{align*}
&(t,x)\mapsto a(t,x,u(t,x),\cR_a u(t,x))\in L^\infty([0,1];C^\infty_b),\\
&(t,x)\mapsto b(t,x,u(t,x),\cR_b u(t,x))\in L^\infty([0,1];C^\infty_b),
\end{align*}
and by (\ref{EO1}),
$$
(t,x)\mapsto f(t,x,u(t,x),\cR_f u(t,x))\in\mY^\infty.
$$
Set $u_0(t,x)\equiv0$. By Theorem \ref{Th1}, we can recursively define $u_n\in\mX^\infty$
by the following linear equation:
\begin{align}
\p_t u_n+a(u_{n-1},\cR_a u_{n-1})(-\Delta)^{\frac{1}{2}} u_n+
b(u_{n-1},\cR_b u_{n-1})\cdot\nabla u_n=f(u_{n-1},\cR_f u_{n-1})\label{EQ2}
\end{align}
subject to the initial value $u_n(0)=\varphi\in\mW^\infty$.

We first assume that $\gamma_0=0$. For simplicity of notation, we set
\begin{align*}
a_n(t,x)&:=a(t,x,u_{n-1}(t,x),\cR_a u_{n-1}(t,x)),\\
b_n(t,x)&:=b(t,x,u_{n-1}(t,x),\cR_b u_{n-1}(t,x)),\\
f_n(t,x)&:=f(t,x,u_{n-1}(t,x),\cR_f u_{n-1}(t,x)).
\end{align*}
By the maximum principle (see Theorem \ref{Max}) and in view of $\gamma_0=1$, it is easy to see that
\begin{align*}
\|u_n(t)\|_\infty&\leq\|\tilde u_n(t)\|_\infty+\int^t_0\|f(s,\cdot,u_{n-1}(s,\cdot), \cR_f u_{n-1}(s,\cdot))\|_\infty\dif s\\
&\leq\|\tilde u_n(0)\|_\infty+\int^t_0(C_{f,0}\|u_{n-1}(s)\|_\infty+\|h_0\|_\infty)\dif s\\
&\leq\|\varphi\|_\infty+\|h_0\|_\infty+C_{f,0}\int^t_0\|u_{n-1}(s)\|_\infty\dif s,
\end{align*}
which yields by Gronwall's inequality that
\begin{align}
\sup_{t\in[0,1]}\|u_n(t)\|_\infty\leq e^{C_{f,0}}(\|\varphi\|_\infty+\|h_0\|_\infty)=:K_0.\label{EU3}
\end{align}
By Theorem \ref{Th6} and (\ref{Ho}), there exist $\beta\in(0,1)$ and constant $K_1>0$ depending on $K_0$
such that for all $n\in\mN$,
$$
\sup_{t\in[0,1]}|u_n(t)|_{\sC^\beta}+\sup_{t\in[0,1]}|\cR_a u_n(t)|_{\sC^\beta}+\sup_{t\in[0,1]}|\cR_b u_n(t)|_{\sC^\beta}\leq K_1.
$$
Thus, by (\ref{EO33}) we have
\begin{align}
&|a_n(t,x)-a_n(t,y)|\leq \|\nabla_x a\|_{L^\infty_{K_0}}|x-y|+K_1\Big(\|\nabla_u a\|_{L^\infty_{K_0}}
+\|\nabla_u a\|_{L^\infty_{K_0}}\Big)|x-y|^\beta,\label{EO4}\\
&|b_n(t,x)-b_n(t,y)|\leq \|\nabla_x b\|_{L^\infty_{K_0}}|x-y|+K_1\Big(\|\nabla_u b\|_{L^\infty_{K_0}}
+\|\nabla_u b\|_{L^\infty_{K_0}}\Big)|x-y|^\beta,\label{EO5}
\end{align}
where $\|\cdot\|_{L^\infty_{K_0}}$ denotes the sup-norm in $L^\infty([0,1]\times\mR^d\times B^m_{K_0}\times\mR^k)$,

For $k=0,1,2,\cdots,$ set
$$
w^{(k)}_n(t,x):=\nabla^k u_n(t,x).
$$
By the chain rule, we have
$$
\p_t w^{(k)}_n+a_n(-\Delta)^{\frac{1}{2}} w^{(k)}_n+b_n\cdot\nabla w^{(k)}_n=g^{(k)}_n,
$$
where $g^{(0)}_n=f_n$ and for $k\geq 1$,
$$
g^{(k)}_n:=\nabla^kf_n-\sum_{j=1}^k \frac{k!}{(k-j)!j!}
\Big(\nabla^ja_n\cdot\nabla^{k-j}(-\Delta)^{\frac{1}{2}} u_n+\nabla^jb_n\cdot\nabla^{k-j}\nabla u_n\Big).
$$
By (\ref{EO4}), (\ref{EO5}) and Lemma \ref{Le7}, we have for all $p>1$ and $t\in[0,1]$,
\begin{align}
\|w^{(k)}_n\|_{\mX^{1,p}_t}\leq C_{k,p}\left(\|\nabla^k\varphi\|_{1-\frac{1}{p},p}+\|g^{(k)}_n\|_{\mY^{0,p}_t}\right),\label{EO6}
\end{align}
where $C_{k,p}$ is independent of $n$.

For $k=0$, by (\ref{EO6}), (\ref{EO1}) and (\ref{EU3}), we have
\begin{align*}
&\|u_n(t)\|^p_p+\int^t_0\|u_n(s)\|^p_{1,p}\dif s\leq C\|\varphi\|^p_{1-\frac{1}{p},p}+C\int^t_0\|f_n(s)\|_p^p\dif s\\
&\qquad\leq C\|\varphi\|^p_{1-\frac{1}{p},p}+\int^t_0\Big(C_{f,0}(K_0^{\gamma_0}+1)
\|u_{n-1}(s)\|_p+\|h_0(s)\|_p\Big)^p\dif s\\
&\qquad\leq C\|\varphi\|^p_{1-\frac{1}{p},p}+C\int^t_0\|u_{n-1}(s)\|^p_p\dif s+C\int^t_0\|h_0(s)\|^p_p\dif s.
\end{align*}
By Gronwall's inequality, one gets
$$
\sup_{n\in\mN}\sup_{t\in[0,1]}\|u_n(t)\|^p_p\leq C_p,
$$
and therefore, for all $p>1$,
$$
\sup_{n\in\mN}\|u_n\|_{\mX^{1,p}}\leq C_p.
$$
Now for any $k=1,2,\cdots$, since by the chain rules, $g^{(k)}_n$ only contains
the powers of all derivatives up to $k$-order of $u_n$, $\cR_a u_n$, $\cR_b u_n$ and $\cR_f u_n$,
by induction method and using H\"older's inequality, it is easy to see
that for all $k\in\mN$ and $p>1$,
\begin{align}
\sup_{n\in\mN}\|u_n\|_{\mX^{k,p}}\leq C_{k,p}.\label{EW2}
\end{align}

Below we write
$$
w_{n,m}(t,x):=u_n(t,x)-u_m(t,x).
$$
Then
$$
\p_t w_{n,m}+a_n(-\Delta)^{\frac{1}{2}} w_{n,m}+ b_n\cdot\nabla w_{n,m}=g_{n,m},
$$
where
$$
g_{n,m}:=f_n-f_m+(a_m-a_n)(-\Delta)^{\frac{1}{2}}u_m+(b_m-b_n)\cdot\nabla u_m.
$$
By Lemma \ref{Le7} again, we have for all $p>1$ and $t\in[0,1]$,
$$
\|w_{n,m}\|_{\mX^{1,p}_t}\leq C\|g_{n,m}\|_{\mY^{0,p}_t}.
$$
Here and below, $C>0$ is independent of $n,m$. Using (\ref{EW2}) and (\ref{EO33}), we have
\begin{align*}
\|g_{n,m}\|_{\mY^{0,p}_t}&\leq C\Big(\|f_n-f_m\|_{\mY^{0,p}_t}+
\|a_n-a_m\|_{\mY^{0,p}_t}+\|b_n-b_m\|_{\mY^{0,p}_t}\Big)\\
&\leq C\Big(\|\nabla_u f\|_{L^\infty_{K_0}}+\|\nabla_r f\|_{L^\infty_{K_0}}+
\|\nabla_u a\|_{L^\infty_{K_0}}+\|\nabla_r a\|_{L^\infty_{K_0}}\\
&\qquad+\|\nabla_u b\|_{L^\infty_{K_0}}+\|\nabla_r b\|_{L^\infty_{K_0}}\Big)\|w_{n-1,m-1}\|_{\mY^{0,p}_t}.
\end{align*}
Hence,
$$
\sup_{s\in[0,t]}\|w_{n,m}(s)\|^p_p\leq C\int^t_0\|w_{n-1,m-1}(s)\|^p_p\dif s.
$$
Taking sup-limits and by Fatou's lemma, we obtain
$$
\varlimsup_{n,m\to\infty}\sup_{s\in[0,t]}\|w_{n,m}(s)\|^p_p\leq C\int^t_0
\varlimsup_{n,m\to\infty}\sup_{s\in[0,r]}\|w_{n,m}(s)\|^p_p\dif r.
$$
So,
$$
\varlimsup_{n,m\to\infty}\sup_{s\in[0,1]}\|w_{n,m}(s)\|^p_p=0,
$$
which together with (\ref{EW2}) and the interpolation inequality yields that for all $k\in\mN$ and $p>1$,
$$
\varlimsup_{n,m\to\infty}\sup_{s\in[0,1]}\|w_{n,m}(s)\|^p_{k,p}=0.
$$
Thus, there exists a $u\in\mX^\infty$ such that for all $k\in\mN$ and $p>1$,
$$
\varlimsup_{n,m\to\infty}\sup_{s\in[0,1]}\|u_n(s)-u(s)\|^p_{k,p}=0.
$$
Taking limits for (\ref{EQ2}), one sees that $u$ solves equation (\ref{EQ1}).

Now we want to drop $\gamma_0=0$ and assume (\ref{EO2}). For $R>0$, let $\chi_R\in C^\infty_0(\mR^d)$
be a nonnegative cutoff function with $\chi_R(u)=1$ for $|u|\leq R$ and $\chi_R(u)=0$ for $|u|>R+1$.
Set
$$
f_R(t,x,u,r):=f(t,x,u,r)\chi_R(u)
$$
Let $u_R\in\mX^\infty$ solve
$$
\p_t u_R+a(u_R,\cR_a u_R)(-\Delta)^{\frac{1}{2}} u_R+ b(u_R,\cR_b u_R)\cdot\nabla u_R=f_R(u_R,\cR_f u_R).
$$
Noticing that by (\ref{For}),
\begin{align*}
2\<(-\Delta)^{\frac{1}{2}} u_R,u_R\>_{\mR^m}=(-\Delta)^{\frac{1}{2}}|u_R|^2+\sE(u_R,u_R),
\end{align*}
we have
\begin{align*}
&2\p_t |u_R|^2+a(u_R,\cR_a u_R)(-\Delta)^{\frac{1}{2}}|u_R|^2+b(u_R,\cR_b u_R)\cdot\nabla|u_R|^2\\
&\quad=2\<u_R,f_R(u_R,\cR_f u_R)\>_{\mR^m}-a(u_R,\cR_a u_R)\sE(u_R,u_R)\stackrel{(\ref{EO2})}{\leq} 2C_f(|u_R|^2+1).
\end{align*}
Thus, by the maximal principle, we have
$$
\|u_R(t)\|^2_\infty\leq\|\varphi\|^2_\infty+C_f\int^t_0(\|u_R(s)\|^2_\infty+1)\dif s,
$$
which implies that for all $R>0$,
$$
\sup_{t\in[0,1]}\|u_R(t)\|^2_\infty\leq e^{C_f}(\|\varphi\|^2_\infty+C_f).
$$
The proof is finished by taking $R:=[e^{C_f}(\|\varphi\|^2_\infty+C_f)]^{1/2}$.
\end{proof}

\section{Fully nonlinear and nonlocal equation: Proof of Theorem \ref{Main}}

The following lemma will play a key role in proving the existence.
\bl\label{Le4}
Let $a\in L^\infty([0,1]; C^1_b(\mR^d))$ be bounded below by $a_0>0$ and $b\in L^\infty([0,1]; C^1_b(\mR^d))$.
Let $u:[0,1]\times\mR^d\to\mR^d$ belong to $\mX^{2,p}$ for some $p>1$ and satisfy
\begin{align}
\p_tu=a(-\Delta)^{-\frac{1}{2}}\square u+b\cdot(\nabla u-(\nabla u)^t),\label{EU4}
\end{align}
where $\square:=\div\nabla-\nabla\div$. Then we have
\begin{align*}
\|u\|_{\mX^{1,p}}+\|U\|_{\mX^{1,p}}\leq C\Big(\|u(0)\|_{1,p}+\|U(0)\|_{1,p}\Big),
\end{align*}
where $U:=\nabla u-(\nabla u)^t$.
\el
\begin{proof}
By equation (\ref{EU4}), one sees that
$$
\p_tu=-a(-\Delta)^{\frac{1}{2}} u-a(-\Delta)^{-\frac{1}{2}}\nabla\div u+b\cdot U,
$$
and
$$
\p_t U=-a(-\Delta)^{\frac{1}{2}}U+b\cdot\nabla U+(\nabla b)\cdot U-[(\nabla b)\cdot U]^t+A,
$$
where
$$
A:=(\nabla a)^t\cdot (-\Delta)^{-\frac{1}{2}}\square u-((-\Delta)^{-\frac{1}{2}}\square u)^t\cdot \nabla a.
$$
By Lemma \ref{Le7}, there exists a constant $C>0$ such that for all $t\in[0,1]$,
\begin{align*}
\|u\|_{\mX^{1,p}_t}&\leq C\|u(0)\|_{1-\frac{1}{p},p}+C\|a(-\Delta)^{-\frac{1}{2}}\nabla\div u\|_{\mY^{0,p}_t}
+C\|b\cdot U\|_{\mY^{0,p}_t}\\
&\leq C\|u(0)\|_{1-\frac{1}{p},p}+C\|a\|_\infty\|\div u\|_{\mY^{0,p}_t}
+C\|b\|_\infty\|U\|_{\mY^{0,p}_t},
\end{align*}
and
\begin{align*}
\|U\|_{\mX^{1,p}_t}&\leq C\|U(0)\|_{1-\frac{1}{p},p}
+C\|(\nabla b)\cdot U+U\cdot(\nabla b)^t+A\|_{\mY^{0,p}_t}\\
&\leq C\|U(0)\|_{1-\frac{1}{p},p}+C(\|\nabla a\|_\infty+\|\nabla b\|_\infty)\|\nabla u\|_{\mY^{0,p}_t}.
\end{align*}
In particular, for all $t\in[0,1]$,
\begin{align}
\|u(t)\|^p_p+\int^t_0\|u(s)\|_{1,p}^p\dif s&\leq
C\|u(0)\|_{1,p}^p+C\int^t_0\|\div u(s)\|^p_p\dif s
+C\int^t_0\|U(s)\|^p_p\dif s\no\\
&\leq C\|u(0)\|_{1,p}^p+Ct\left(\sup_{s\in[0,t]}\|\div u(s)\|^p_p
+\sup_{s\in[0,t]}\|U(s)\|^p_p\right),\label{EU2}
\end{align}
and
\begin{align}
\|U(t)\|^p_p+\int^t_0\|U(s)\|^p_{1,p}\dif s\leq
C\|U(0)\|_{1,p}^p+C\int^t_0\|u(s)\|^p_{1,p}\dif s.\label{EU7}
\end{align}
On the other hand, noticing that
$$
\div\square u=0,
$$
we have
$$
\p_t\div u=\nabla a\cdot(-\Delta)^{-\frac{1}{2}}\square u+\div(b\cdot U).
$$
Hence,
\begin{align}
\|\div u(t)\|_p&\leq \|\div u(0)\|_p+\|\nabla a\|_\infty\int^t_0\|(-\Delta)^{-\frac{1}{2}}\square u(s)\|_p\dif s\no\\
&\quad+\|b\|_\infty\int^t_0\|\nabla U(s)\|_p\dif s+\|\nabla b\|_\infty\int^t_0\|U(s)\|_p\dif s\no\\
&\stackrel{(\ref{EU7})}{\leq}
C\Big(\|\div u(0)\|_p+\|U(0)\|_{1,p}\Big)+C\left(\int^t_0\|u(s)\|^p_{1,p}\dif s\right)^{1/p}.\label{EU6}
\end{align}
Now substituting (\ref{EU7}) and (\ref{EU6}) into (\ref{EU2}), we obtain that for all $t\in[0,1]$,
$$
\|u(t)\|_p^p+\int^t_0\|u(s)\|_{1,p}^p\dif s\leq C_0\Big(\|u(0)\|^p_{1,p}+\|U(0)\|^p_{1,p}\Big)+
C_1 t\int^t_0\|u(s)\|^p_{1,p}\dif s,
$$
where $C_0, C_1$ are independent of $\|u(0)\|_{1,p}$ and $\|U(0)\|_{1,p}$.
Choosing $t_0:=1/(2C_1)$, we arrive at
$$
\sup_{t\in[0,t_0]}\|u(t)\|_p^p+\int^{t_0}_0\|u(s)\|_{1,p}^p\dif s\leq 2C_0\Big(\|u(0)\|^p_{1,p}+\|U(0)\|^p_{1,p}\Big).
$$
So, for some $C_2>0$,
$$
\|u\|^p_{\mX^{1,p}_{t_0}}+\|U\|^p_{\mX^{1,p}_{t_0}}\leq C_2\Big(\|u(0)\|^p_{1,p}+\|U(0)\|^p_{1,p}\Big).
$$
In particular,
$$
\int^{t_0}_{2t_0/3}\Big(\|u(s)\|_{1,p}^p+\|U(s)\|_{1,p}^p\Big)\dif s\leq C_2\Big(\|u(0)\|^p_{1,p}+\|U(0)\|^p_{1,p}\Big).
$$
Thus, there is at least one point $s_0\in[2t_0/3,t_0]$ such that
$$
\|u(s_0)\|_{1,p}^p+\|U(s_0)\|_{1,p}^p\leq \frac{3C_2}{t_0}\Big(\|u(0)\|^p_{1,p}+\|U(0)\|^p_{1,p}\Big).
$$
Now starting from $s_0$, as above, one can prove that for the same $t_0$,
$$
\|u(\cdot+s_0)\|^p_{\mX^{1,p}_{t_0}}+\|U(\cdot+s_0)\|^p_{\mX^{1,p}_{t_0}}
\leq C_2\Big(\|u(s_0)\|^p_{1,p}+\|U(s_0)\|^p_{1,p}\Big)
\leq \frac{3C_2^2}{t_0}\Big(\|u(0)\|^p_{1,p}+\|U(0)\|^p_{1,p}\Big).
$$
Repeating the above proof, we obtain the desired estimate.
\end{proof}

We are now in a position to give

{\bf Proof of Theorem \ref{Main}}: We divide the proof into three steps.

{\bf (Step 1)}. In this step we consider the following fully non-linear and nonlocal parabolic equation:
$$
\p_t u=F(t,x,\nabla u,-(-\Delta)^{\frac{1}{2}} u),\ \ u(0)=\varphi.
$$
As introduced in the introduction, let
\begin{align}
\cR w=(-\Delta)^{-\frac{1}{2}}\div w.\label{LP1}
\end{align}
For any $\varphi\in\mU^\infty=\cap_{k,p}\mU^{k,p}$, where $\mU^{k,p}$ is defined by (\ref{EM7}),
by Theorem \ref{Th2}, there exists a unique $w\in\mX^\infty$ solving the following parabolic system:
$$
\p_tw=-(\p_qF)(w,\cR w)(-\Delta)^{\frac{1}{2}}w+(\nabla_wF)(w,\cR w)\nabla w+\nabla_x F(w,\cR w)
$$
subject to $w(0)=\nabla\varphi$. Define
$$
u(t,x):=\varphi(x)+\int^t_0F(s,x,w(s,x),\cR w(s,x))\dif s
$$
and
$$
h(t,x):=\nabla u(t,x)-w(t,x).
$$
Then we have
\begin{align*}
\p_t h&=(\p_qF)(w,\cR w)(\nabla\cR w+(-\Delta)^{\frac{1}{2}}w)+(\nabla_w F)(w,\cR w)((\nabla w)^t-\nabla w)\\
&=(\p_qF)(w,\cR w)(-\Delta)^{-\frac{1}{2}}\square h
+(\nabla_w F)(w,\cR w)(\nabla h-(\nabla h)^t)
\end{align*}
subject to $h(0)=0$, where $\square:=\div\nabla-\nabla\div$.
By Lemma \ref{Le4}, we have
$$
h=0\Rightarrow w=\nabla u.
$$
Thus, by (\ref{LP1}),
$$
\p_tu(t,x)=F(t,x,\nabla u(t,x),\cR \nabla u(t,x))
=F(t,x,\nabla u(t,x),-(-\Delta)^{\frac{1}{2}}u(t,x)).
$$
By the maximum principle (see Theorem \ref{Max}), we have
\begin{align}
\|u(t)\|_\infty\leq\|\varphi\|_\infty+\int^t_0\|F(s,\cdot,0,0)\|_\infty\dif s.\label{ED1}
\end{align}
In particular, $u\in C([0,1];\mU^\infty)$.

{\bf (Step 2)}.  Now we consider the general case.
Set $u_0(t,x)=0$. Let $u_n\in C([0,1];\mU^\infty)$ be defined recursively
by the following equation:
\begin{align}
\p_tu_n=F(t,x,u_{n-1},\nabla u_n,-(-\Delta)^{\frac{1}{2}} u_n),\ \ u_n(0)=\varphi.\label{EH5}
\end{align}
By (\ref{ED1}) and (\ref{EW3}), we have
\begin{align*}
\|u_n(t)\|_\infty&\leq\|\varphi\|_\infty+\int^t_0\|F(s,\cdot,u_{n-1}(s,\cdot),0,0)\|_\infty\dif s\\
&\leq\|\varphi\|_\infty+\kappa_0\int^t_0(\|u_{n-1}(s)\|_\infty+1)\dif s.
\end{align*}
By Gronwall's inequality, we get
\begin{align}
\|u_n(t)\|_\infty\leq e^{\kappa_0}(\|\varphi\|_\infty+\kappa_0)=:K_0.\label{EW10}
\end{align}
On the other hand, by taking gradients with respect to $x$ for equation (\ref{EH5}), we have
\begin{align*}
\p_t\nabla u_n&=-\p_q F(t,x,u_{n-1},\nabla u_n,-(-\Delta)^{\frac{1}{2}} u_n)(-\Delta)^{\frac{1}{2}}\nabla u_n\\
&\quad+\nabla_w F(t,x,u_{n-1},\nabla u_n,-(-\Delta)^{\frac{1}{2}} u_n)\nabla^2 u_n\\
&\quad+\p_u F(t,x,u_{n-1},\nabla u_n,-(-\Delta)^{\frac{1}{2}} u_n)\nabla u_{n-1}\\
&\quad+\nabla_x F(t,x,u_{n-1},\nabla u_n,-(-\Delta)^{\frac{1}{2}} u_n).
\end{align*}
By the maximum principle again and (\ref{EW6}), (\ref{EW7}) with $\gamma_{K_0,1}=0$,  we have
\begin{align*}
\|\nabla u_n(t)\|_\infty&\leq\|\nabla\varphi\|_\infty+\int^t_0\|\p_u F(s,x,u_{n-1},
\nabla u_n,-(-\Delta)^{\frac{1}{2}} u_n)\nabla u_{n-1}\|_\infty\dif s\\
&\quad+\int^t_0\|\nabla_x F(s,x,u_{n-1},\nabla u_n,
-(-\Delta)^{\frac{1}{2}} u_n)\|_\infty\dif s\\
&\leq\|\nabla\varphi\|_\infty+C\int^t_0\Big(\|\nabla u_{n-1}(s)\|_\infty+\|\nabla u_n(s)\|_\infty+1\Big)\dif s,
\end{align*}
where $C$ is independent of $n$. By Gronwall's inequality, we get
\begin{align}
\sup_n\sup_{t\in[0,1]}\|\nabla u_n(t)\|_\infty<+\infty.\label{EH1}
\end{align}
Moreover, by (\ref{EW5}), (\ref{EW6}), (\ref{EW7}), (\ref{EW10}), Theorem \ref{Th6} and Lemma \ref{Le1},
as in the proof of Theorem \ref{Th2},
we have for all $p>1$,
$$
\|\nabla u_n\|_{\mX^{1,p}_t}\leq C\left(\|\nabla\varphi\|_{1-\frac{1}{p},p}+\|\nabla u_{n-1}\|_{\mY^{0,p}_t}
+\|\nabla u_n\|_{\mY^{0,p}_t}+\|h_1\|_p\right),
$$
which implies by Gronwall's inequality that
\begin{align}
\sup_n\|\nabla u_n\|_{\mX^{1,p}}<+\infty,\label{EH2}
\end{align}
and furthermore, for all $k\in\mN$ and $p>1$,
\begin{align}
\sup_n\|\nabla u_n\|_{\mX^{k,p}}<+\infty.\label{EH4}
\end{align}
This together with (\ref{EW10}) gives
\begin{align}
\sup_n\sup_{t\in[0,1]}\|u_n(t)\|_{\mU^{k,p}}<+\infty.\label{EH44}
\end{align}

{\bf (Step 3)}. Next we want to show that $u_n$ converges to some $u$ in $C([0,1];\mU^{k,p})$. For $n,m\in\mN$, set
$$
v_{n,m}(t,x):=u_n(t,x)-u_m(t,x).
$$
Then
$$
\p_t v_{n,m}=-a_{n,m}(-\Delta)^{\frac{1}{2}}v_{n,m}+b_{n,m}\cdot\nabla v_{n,m}+f_{n,m}v_{n-1,m-1},
$$
where
\begin{align*}
a_{n,m}&:=\int^1_0(\p_qF)(u_{n-1},\nabla u_n,-(-\Delta)^{\frac{1}{2}}(sv_{n,m}+u_m))\dif s,\\
b_{n,m}&:=\int^1_0(\nabla_wF)(u_{n-1},\nabla (sv_{n,m}+u_m),-(-\Delta)^{\frac{1}{2}}u_m)\dif s,\\
f_{n,m}&:=\int^1_0(\p_uF)(sv_{n-1,m-1}+u_{m-1},\nabla u_m,-(-\Delta)^{\frac{1}{2}}u_m)\dif s.
\end{align*}
By the maximum principle, we have
$$
\|v_{n,m}(t)\|_\infty\leq C\int^t_0\|v_{n-1,m-1}(s)\|_\infty\dif s,
$$
and by Gronwall's inequality,
\begin{align}
\lim_{n,m\to\infty}\sup_{t\in[0,1]}\|v_{n,m}(t)\|_\infty=0.\label{EK10}
\end{align}
On the other hand, by Lemma \ref{Le1} and (\ref{EH44}), we may derive that for all $t\in[0,1]$,
$$
\|v_{n,m}\|_{\mX^{1,p}_t}\leq C\|v_{n-1,m-1}\|_{\mY^{0,p}_t},
$$
and so,
$$
\lim_{n,m\to\infty}\|v_{n,m}\|_{\mX^{1,p}}=0.
$$
This together with (\ref{EH4}), the interpolation inequality and (\ref{EK10}) yields that for all $k\in\mN$ and $p>1$,
$$
\lim_{n,m\to\infty}\sup_{t\in[0,1]}\|v_{n,m}(t)\|_{\mU^{k,p}}=0.
$$
Thus, there is a $u\in C([0,1];\mU^\infty)$ such that  for all $k\in\mN$ and $p>1$,
$$
\lim_{n,m\to\infty}\sup_{t\in[0,1]}\|u_n(t)-u(t)\|_{\mU^{k,p}}=0.
$$
The proof is complete by taking limits for approximation equation (\ref{EH5}).

\end{document}